\documentclass[leqno]{amsart}
\usepackage{amssymb}
\usepackage{euscript}
\usepackage[hypertex]{hyperref}

\title[Riemannian symmetries in flag manifolds]{Riemannian symmetries in flag manifolds}

\author[Piu - Remm]{Paola Piu and Elisabeth Remm}
\thanks{The first author was supported by a Visiting Professor fellowship at Universit\'e de Haute Alsace - Mulhouse in February 2012  and by GNSAGA(Italy)}
\address{Dipartimento di Matematica e Informatica, Universit\`a degli Studi di Cagliari, Via Ospedale 72, 09124 Cagliari, ITALIA}
\email{piu@unica.it}


\address{Laboratoire de Math\'ematiques et Applications,
        Universit\'e de Haute Alsace, Facult\'e des Sciences et
        Techniques, 4, rue des Fr\`eres Lumi\`ere,
        68093~Mulhouse~cedex, France.}
\email{Elisabeth.Remm@uha.fr}

\begin{document}

\newtheorem{theorem}{Theorem}
\newtheorem{acknowledgement}[theorem]{Acknowledgement}
\newtheorem{algorithm}[theorem]{Algorithm}
\newtheorem{axiom}[theorem]{Axiom}
\newtheorem{case}[theorem]{Case}
\newtheorem{claim}[theorem]{Claim}
\newtheorem{conclusion}[theorem]{Conclusion}
\newtheorem{condition}[theorem]{Condition}
\newtheorem{conjecture}[theorem]{Conjecture}
\newtheorem{corollary}[theorem]{Corollary}
\newtheorem{criterion}[theorem]{Criterion}
\newtheorem{definition}[theorem]{Definition}
\newtheorem{example}[theorem]{Example}
\newtheorem{exercise}[theorem]{Exercise}
\newtheorem{lemma}[theorem]{Lemma}
\newtheorem{notation}[theorem]{Notation}
\newtheorem{problem}[theorem]{Problem}
\newtheorem{proposition}[theorem]{Proposition}
\newtheorem{remark}[theorem]{Remark}
\newtheorem{solution}[theorem]{Solution}
\newtheorem{summary}[theorem]{Summary}
\newcommand\R{\mathbb{R}}
\newcommand\g{\frak{g}}
\newcommand\K{\mathbb{K}}
\newcommand\C{\mathbb{C}}
\newcommand\Z{\mathbb{Z}}
\newcommand{\zz}{$\Z^2_2$-symmetric space}
\newcommand\z{ \mathbb{Z}_2^2}
\newcommand\p{\mathcal{P}}
\newcommand{\pf}{{\it Proof. }}
\newcommand\f{\frak{t}}
\newcommand\h{\frak{h}}
\newcommand\m{\frak{m}}
\setlength{\oddsidemargin}{0in}
\setlength{\evensidemargin}{.25in}
\setlength{\textwidth}{6.25in}

\pagestyle{myheadings}

\begin{abstract}
Flag manifolds are  in general not symmetric spaces. But they are  provided with a structure of $\mathbb{Z}_2^k$-symmetric space. We describe the Riemannian metrics adapted to this structure and some properties of reducibility. We detail for the flag manifold $SO(5)/SO(2)\times SO(2) \times SO(1)$ what are the conditions for a metric adapted to the $\mathbb{Z}_2^2$-symmetric structure to be naturally reductive.
\end{abstract}

\maketitle

\bibliographystyle{plain}

\section{Introduction}

In this work we call flag manifold any homogeneous space
$$G/H=\frac{SO(n)}{SO(p_1) \times \cdots \times SO(p_k) }$$
with $p_1 + \cdots +p_k=n.$ We suppose also that $p_1\geq p_2 \geq \cdots\geq p_k.$
When $k=2$ and $p_2=1,$ this space is isomorphic to the sphere $S^{n-1}$ and it is a symmetric space. When $k=2$ and $p_1 \neq 1$, the homogeneous space $G/H$ is isomorphic to the Grassmannian
manifold and it is also a symmetric space. But if $k>2$ the homogeneous space $G/H$ is reductive but not symmetric. In \cite{Ba-Go} and \cite{Go-Re-Santiago} we have shown that we can define on $G/H$ a structure of $\mathbb{Z}_2^k$-symmetric space, that is, the Lie algebra $\frak{g}$ of $G$ admits a $\mathbb{Z}_2^k$-grading. A Riemannian metric on $G/H$ is called associated with the $\mathbb{Z}_2^k$-symmetric structure if the natural symmetries defining the $\mathbb{Z}_2^k$-symmetric structure are isometries. The Riemannian geometry is more complicated than the Riemannian geometry of symmetric spaces (or $\mathbb{Z}_2$-symmetric spaces); in particular, it is not always naturally reductive. In the paper we investigate the $\mathbb{Z}_2^k$-symmetric Riemannian tensor  on flag manifolds and develop all the computations when $G/H=SO(5)/SO(2)\times SO(2) \times SO(1).$

\section{Riemannian
$\mathbb{Z}_2^k$-symmetric spaces}
This notion has been introduced in \cite{Lutz}. Let $G/H$ be an homogeneous space with a connected Lie group $G$. We denote by $\frak{g}$ and $\frak{h}$ respectively the Lie algebras of $G$ and $H$.

\begin{definition}
The homogeneous space $G/H$ is $\mathbb{Z}_2^k$-symmetric if the Lie algebra $\frak{g}$ admits a $\mathbb{Z}_2^k$-grading, that is,
$$\frak{g}=\oplus_{\gamma \in \mathbb{Z}_2^k}\frak{g}_{\gamma}, \quad [\frak{g}_{\gamma_1},\frak{g}_{\gamma_2}]\subset \frak{g}_{\gamma_1\gamma_2},$$
with $\frak{g}_{e}=\frak{h}$ where $e$ is the identity of $\mathbb{Z}_2^k$.
\end{definition}

If $\frak{m}=\oplus_{\gamma \in \mathbb{Z}_2^k, \gamma \neq e} \frak{g}_{\gamma}$ then $\frak{g}=\frak{g}_{e} \oplus \frak{m}$ and the $\mathbb{Z}_2^k$-grading implies $[\frak{g}_e, \frak{m}]\subset \frak{m}.$
But if $k \geq 2$, we have not $[\frak{m},\frak{m}]\subset \frak{g}_e.$
Thus the decomposition $\frak{g}= \frak{g}_e \oplus \frak{m}$ is reductive, not symmetric when $k\geq 2$.

\smallskip

\noindent {\bf Consequence.} A $\mathbb{Z}_2^k$-symmetric homogeneous space $G/H$ is reductive.

\begin{example} From \cite{Ba-Go}, \cite{Go-Re-Santiago} and \cite{Kollross},  it is possible to give a classification of the  $\mathbb{Z}_2^2$-grading of complex simple Lie algebras. In the following list we give the pairs $(\frak{g}, \frak{g}_e=\frak{h})$
which the (local) classification of $\mathbb{Z}_2^2$-symmetric structures when $G$ is a simple complex or $G$ is simple compact and real.
$$
\begin{array}{l|l}
\frak{g} & \frak{g}_e=\frak{h}  \\
\hline \\
so(k_1+k_2+k_3+k_4),  & \oplus_{i=1}^4 so(k_i) \\
k_1 \geq k_2 \geq k_3 \geq k_4, k_3 \neq 0   \\
sp(k_1+k_2+k_3+k_4),  &  \oplus_{i=1}^4 sp(k_i)  \\
k_1 \geq k_2 \geq k_3 \geq k_4, k_3 \neq 0 \\
so(2m) & gl(m) \\
so(2k_1+2k_2) & gl(k_1)\oplus gl(k_2) \\
sp(2m) & gl(m) \\
sp(2k_1+2k_2) & gl(k_1)\oplus gl(k_2) \\
so(2m) & so(m)\\
so(4m) & sp(2m) \\
sp(4m) & sp(2m) \\
sp(2m) & so(m) \\
sl(2n) & sl(n) \\
sl(k_1+k_2) & \oplus_{i=1}^2 sl(k_i) \oplus \mathbb{C} \\
sl(k_1+k_2+k_3) & \oplus_{i=1}^3 sl(k_i)\oplus \mathbb{C}^2 \\
sl(k_1+k_2+k_3+k_4) & \oplus_{i=1}^4 sl(k_i) \oplus   \mathbb{C}^3 \\
so(8) & gl(3) \oplus gl(1)
\end{array}
\quad
\begin{array}{l|l}
\frak{g} & \frak{g}_e=\frak{h}  \\
\hline \\
E_6 & so(6)\oplus \mathbb{C}\\
    & sp(2)\oplus sp(2) \\
    & sp(3)\oplus sp(1) \\
    & su(3)\oplus su(3)\oplus \mathbb{C}^2 \\
    & su(4)\oplus sp(1)\oplus sp(1)\oplus \mathbb{C} \\
    & su(5) \oplus \mathbb{C}^2 \\
    & so(8) \oplus \mathbb{C}^2 \\
    & so(9) \\
    \hline \\
E_7 & so(8) \\
    & su(4)\oplus su(4)\oplus \mathbb{C} \\
    & sp(4) \\
    & su(6)\oplus sp(1)\oplus \mathbb{C} \\
    & so(8)\oplus so(4)\oplus sp(1) \\
    & u(6) \oplus \mathbb{C} \\
    & so(10) \oplus \mathbb{C}^2 \\
    & F_4
     \end{array}
 $$
 $$
 \begin{array}{l|l}
\frak{g} & \frak{g}_e=\frak{h}  \\
\hline \\
E_8 & so(8) \oplus so(8) \\
    & su(8)\oplus \mathbb{C} \\
    & so(12)\oplus sp(1)\oplus sp(1) \\
    & E_6 \oplus \mathbb{C}^2
    \end{array}
    \quad
\begin{array}{l|l}
\frak{g} & \frak{g}_e=\frak{h}  \\
\hline \\
F_4 & u(3) \oplus \mathbb{C} \\
    & sp(2)\oplus sp(1)\oplus sp(1) \\
    & so(8) \\
G_2 & \mathbb{C}^2
\end{array}
$$
\end{example}

\medskip

If $M=G/H$ is a $\mathbb{Z}_2^k$-symmetric space, the grading
$$\frak{m}=\oplus_{\gamma \in \mathbb{Z}_2^k, \gamma \neq e} \frak{g}_{\gamma}$$
is associated with a spectral decomposition of $\mathfrak{g}$ defined by a family
$\left\{ \sigma_\gamma, \gamma \in \mathbb{Z}_2^k \right\}$ of automorphisms of $\frak{g}$ satisfying
$$\left\{
\begin{array}{l}
\sigma_{\gamma}^2=Id, \\
\sigma_{\gamma_1} \circ \sigma_{\gamma_2}=\sigma_{\gamma_2} \circ \sigma_{\gamma_1},
\end{array}
\right.
$$
for any $\sigma,\sigma_1,\sigma_2 \in \mathbb{Z}_2^k.$ Any $\sigma_{\gamma} \in Aut(\frak{g})$ defines an automorphism $s_\gamma$ of $G$ and $H$, if it is connected, corresponds to the identity component of the group $\left\{A \in G / \forall \sigma \in \mathbb{Z}_2^k,  s_\gamma(A)=A \right\}.$ The family
$\left\{ s_\gamma , \gamma \in \mathbb{Z}_2^k \right\}$ is a subgroup of $Aut(G)$
 isomorphic to $\mathbb{Z}_2^k.$ For any $x \in G/H$, it determines a subgroup $\left\{ s_{\gamma,x} , \gamma \in \mathbb{Z}_2^k \right\}$ of $Diff(M)$ isomorphic to $\mathbb{Z}_2^k.$ The diffeomorphisms $s_{\gamma,x}$ are called the symmetries of the $\mathbb{Z}_2^k$-symmetric space $G/H.$
By extension, we will also call symmetries, the automorphisms $s_\gamma$ of $G.$

\noindent{\bf Remark.} A $\Gamma$-symmetric space, when $\Gamma$ is a cyclic group is usually called generalized symmetric space. In this case the grading of the Lie algebra is defined in the complex field and correspond to the roots of the unity. For a general presentation, see \cite{Kowalski}.
\begin{definition}
A Riemannian metric $g$ on the (reductive) $\mathbb{Z}_2^k$-symmetric space $G/H$ is called adapted to the $\mathbb{Z}_2^k$-symmetric structure if the symmetries $s_{\gamma,x}$ are isometries for any $\gamma$  in $ \mathbb{Z}_2^k$ and $x$ in $G/H.$

In this case we will say that $G/H$ is a Riemannian $\mathbb{Z}_2^k$-symmetric space.
\end{definition}

As $G/H$ is a reductive homogeneous space, we consider only $G$-invariant Riemannian metrics on $G/H$ which are positive-definite $ad(H)$-invariant  symmetric bilinear form $B$ on $\frak{m}$, with the following correspondence
$$B(X,Y)=g(X,Y)_0$$
for $X,Y \in \frak{m}$. As in this paper we consider $H$ connected, the invariance of $B$ is written
$$B([Z,X],Y)+B(X,[Z,Y])=0,$$ for $X,Y$ in $\frak{m}$ and $Z \in \frak{h}.$

\begin{proposition}
Let $G/H$ be a $\mathbb{Z}_2^k$-symmetric structure with $G$ and $H$ connected. Any Riemannian metric $g$ adapted to the $\mathbb{Z}_2^k$-symmetric structure is in one-to-one correspondence with the $ad(H)$-invariant positive-definite bilinear form $B$ on $\frak{m}$ such that $B(\frak{g}_\gamma, \frak{g}_{\gamma '})=0$ for  $\gamma \neq \gamma '$ in $\mathbb{Z}_2^k$ where
$\frak{g}=\oplus_{\gamma \in \mathbb{Z}_2^k}\frak{g}_{\gamma}$ is the $\mathbb{Z}_2^k$-grading corresponding to the $\mathbb{Z}_2^k$-symmetric structure of $G/H$.
Let $U(X,Y)$ be the symmetric bilinear mapping of $\frak{m} \times \frak{m}$ in $\frak{m}$ defined by
$$2B(U(X,Y),Z)=B(X,[Z,Y]_\frak{m})+B([Z,X]_\frak{m},Y),$$
for all $X,Y,Z \in \frak{m},$ where $[\ , ]_\frak{m}$ denote the projection on $\frak{m}$ of the bracket of $\frak{g}.$
\end{proposition}
Thus the Riemmanian connection for $g$ is given by
$$\nabla_X Y=U(X,Y)+\frac{1}{2}[X,Y]_\frak{m}$$
and the curvature tensor satisfies
$$R(X,Y)(Z)=\nabla_X\nabla_Y Z-\nabla _Y \nabla_X Z- \nabla_{[X,Y]_{\m}}Z-[[X,Y]_\frak{\h},Z],$$
for $X,Y,Z \in \m$ and  where the term $\left[\left[X,Y\right]_\frak{\h},Z\right]$ corresponds to the linear isotropy representation of $H$ into $G/H.$
\section{The Riemannian \zz \\
$SO(5)/SO(2)\times SO(2) \times SO(1)$}
Any Riemannian symmetric space is naturally reductive. It is not usually the case for Riemannian $\mathbb{Z}_2^k$-symmetric spaces as soon as $k \geq 2.$ In this section we describe the Riemannian and Ricci tensors for any $\Z^2$-symmetric metric on the flag manifold $SO(5)/SO(2)\times SO(2) \times SO(1)$. In particular we study some properties of these metrics when they are not naturally reductive. A $\Z^2_2$-grading of the Lie algebra $so(5)$ is given by the decomposition
$$\g=\g_e \oplus \g_a \oplus \g_b \oplus \g_c,$$
corresponding to the multiplication of $\Z^2_2$ : $a^2=b^2= c^2=e, \ ab=c, ac=b, bc=a$, $e$ being the identity. To describe the components of the grading, we consider
$$so(5)=\left\{
\left(
\begin{array}{ccccc}
0 & x_1 & a_1 & a_2 & b_1 \\
-x_1 & 0 & a_3 & a_4 & b_2 \\
-a_1 & -a_3 & 0 & x_2 & c_1 \\
-a_2 & -a_4 & -x_2 & 0 & c_2 \\
-b_1 & -b_2 & -c_1 & -c_2 & 0
\end{array}
\right),x_i,a_i,b_i,c_i \in \R\right\}.
$$
Thus $$\begin{array}{l}
\g_e=\left\{ X \in so(5)\,  /  \, a_i=b_i=c_i=0 \right\}, \\
\g_a=\left\{ X \in so(5)\,  /  \, x_i=b_i=c_i=0 \right\}, \\
\g_b=\left\{ X \in so(5)\,  /  \, x_i=b_i=c_i=0 \right\}, \\
\g_c=\left\{ X \in so(5)\,  /  \, x_i=a_i=b_i=0 \right\}.
\end{array}$$
We have $\g_e=so(2) \oplus so(2)\oplus so(1)$ and this grading induces the $\z$-symmetric structure on $SO(5)/SO(2)\times SO(2) \times SO(1)$. Moreover, from \cite{Ba-Go}, this grading is unique up an equivalence of $\z$-gradings.
We denote by $\left\{ \left\{ X_1,X_2 \right\}, \left\{ A_1,A_2,A_3,A_4 \right\},\right.$
$\left. \left\{ B_1,B_2 \right\}, \left\{ C_1,C_2 \right\} \right\}$ the basis of $so(5)$ where each big letter corresponds to
the matrix of $so(5)$ with the small letter equal to $1$ and other  coefficients are zero. This basis is
adapted to the grading. Let us denote by
$\left\{ \omega_1, \omega_2, \alpha_1, \alpha_2, \alpha_3, \alpha_4, \beta_1, \beta_2, \gamma_1, \gamma_2  \right\}$
the dual basis.

\medskip

\begin{theorem}
Any $\Z^2_2$-symmetric Riemannian metric $g$ on $SO(5)/SO(2)\times SO(2) \times SO(1)$ is given by an $ad(H)$-invariant symmetric bilinear form $B$ on $\m=\g_a \oplus \g_b \oplus \g_c$
$$B=t^2(\alpha_1^2+ \alpha_2^2+\alpha_3^2+\alpha_4^2)+u(\alpha_1 \alpha_4-\alpha_2 \alpha_3 )+
v^2(\beta_1^2 +\beta_2^2)+w^2( \gamma_1^2+ \gamma_2^2)$$
with $tvw \ne 0$ and $u \in \left]-4t^2, 4t^2\right[$.
\end{theorem}
\noindent \pf It is explain in detail in \cite{Go-Re-Santiago}.

\medskip

The Riemannian connection is given by
$$\bigtriangledown _X (Y)= \frac{1}{2}[X,Y]_{\m} +U(X,Y),$$
for any $X,Y \in \m,$  where $U$ is the symmetric bilinear mapping on $\m \times \m$ into $\m$ defined by
$$2B(U(X,Y),Z)=B(X,[Z,Y]_{\m})+B([Z,X]_{\m},Y),$$
for any $X,Y,Z  \in \m.$ The bilinear mapping $U$ is reduced to
$$B=\sum_{i=1}^4 \widetilde{\alpha}_i^2+\sum_{i=1}^2 \widetilde{\beta}_i^2+\sum_{i=1}^2 \widetilde{\gamma}_i^2$$
with
$$
\begin{array}{llll}
\widetilde{\alpha}_1=t\alpha_1+\frac{u}{2t}\alpha_4, & \widetilde{\alpha}_2=t\alpha_2-\frac{u}{2t}\alpha_3,& \widetilde{\alpha}_3=K\alpha_3, & \widetilde{\alpha}_4=K\alpha_4, \\
\widetilde{\beta}_1=v\beta_1, &\widetilde{\beta}_2=v\beta_2,&\widetilde{\gamma}_1=w\gamma_1,&\widetilde{\gamma}_2=w\gamma_2.
\end{array}
$$
with $K=\sqrt{t^2-\frac{u^2}{4t^2}}.$ If $\{\widetilde{A_j},\widetilde{B_i},\widetilde{C_i} \}$ is the dual basis, this basis is orthonormal and the brackets $[ \  , \ ]_{\m} $ are given by
$$
\begin{array}{c|c|c|c|c|c|c|c|c}
    &  \widetilde{A_1}  &  \widetilde{A_2} & \widetilde{A_3} &\widetilde{A_4} & \widetilde{B_1} & \widetilde{B_2} & \widetilde{C_1} & \widetilde{C_2} \\
\hline
 & & & & & & & & \\
\widetilde{A_1}     & 0    & 0 &0 & 0   &-\frac{w}{tv}\widetilde{C_1} &  0  & \frac{v}{tw}\widetilde{B_1} & 0   \\
 & & & & & & & & \\
\widetilde{A_2}      &      & 0    &0    & 0&-\frac{w}{tv}\widetilde{C_2} &  0  & 0   & \frac{v}{tw}\widetilde{B_1} \\
 & & & & & & & & \\
\widetilde{A_3}    &      &      & 0   & 0& -\frac{uw}{2t^2vK}\widetilde{C_2}   &-\frac{w}{Kv}\widetilde{C_1} & \frac{v}{Kw}\widetilde{B_2} & \frac{uv}{2t^2wK}\widetilde{B_1}   \\
 & & & & & & & & \\
\widetilde{A_4}      &      &      &    &  0   & \frac{uw}{2t^2vK}\widetilde{C_1}   &-\frac{w}{Kv}\widetilde{C_2}& -\frac{uv}{2t^2wK}\widetilde{B_1}  & \frac{v}{Kw}\widetilde{B_2} \\
 & & & & & & & & \\
\widetilde{B_1}     &      &      &    &      & 0   &0 &-\frac{t}{vw}\widetilde{A_1} &-\frac{t}{vw}\widetilde{A_2} \\
 & & & & & & & & \\
 \widetilde{B_2}    &      &      &    &      &     & 0   &\frac{u}{2vwt}\widetilde{A_2}- &-\frac{u}{2vwt}\widetilde{A_1}- \\
  & & & & & & &\frac{K}{vw}\widetilde{A_3} & \frac{K}{vw}\widetilde{A_4}\\
 & & & & & & & & \\
\widetilde{C_1}     &      &      &    &      &     &     & 0   &0 \\
 & & & & & & & & \\
\widetilde{C_2}   &      &      &    &      &     &     &     &  0  \\
\end{array}
$$
Recall that the Riemannian connection $\nabla$ is given by
$$\nabla_X(Y)=U(X,Y)+\frac{1}{2}[X,Y]_{\m},$$ and the Riemannian curvature $R(X,Y)$ is the matrix given by
$$R(X,Y)(Z)=\nabla_X\nabla_Y Z-\nabla _Y \nabla_X Z- \nabla_{[X,Y]_{\m}}Z-[[X,Y]_\frak{\h},Z],$$
for $X,Y,Z \in \m$ and  where the term $\left[\left[X,Y\right]_\frak{\h},Z\right]$ corresponds to the linear isotropy representation of $H$ into $G/H.$

The symmetric mapping $U$ is given by:
$$\begin{array}{c|c|c|c|c|c|c|c|c}
  U(X,Y)  &  \widetilde{A_1}  &  \widetilde{A_2} & \widetilde{A_3} &\widetilde{A_4} & \widetilde{B_1} & \widetilde{B_2} & \widetilde{C_1} & \widetilde{C_2} \\
\hline
 & & & & & & & & \\
 \widetilde{A_1}     & 0    & 0 &0 & 0   &\frac{t^2-v^2}{2tvw}\widetilde{C_1} &  \frac{u}{4vwt}\widetilde{C_2}  & \frac{-t^2+w^2}{2tvw}\widetilde{B_1} & -\frac{u}{4vwt}\widetilde{B_2}   \\
 & & & & & & & & \\
\widetilde{A_2}      &      & 0    &0    & 0& \frac{t^2-v^2}{2tvw}\widetilde{C_2} &  -\frac{u}{4vwt}\widetilde{C_1}  & \frac{u}{4vwt}\widetilde{B_2}   & \frac{-t^2+w^2}{2tvw}\widetilde{B_1} \\
 & & & & & & & & \\
\widetilde{A_3}    &      &      & 0   & 0& \frac{-uv}{4t^2wK}\widetilde{C_2}    &\frac{K^2-v^2}{2Kvw}\widetilde{C_1} & \frac{-K^2+w^2}{2Kvw}\widetilde{B_2} & \frac{uw}{4t^2vK}\widetilde{B_1}  \\
 & & & & & & & & \\
\widetilde{A_4}      &      &      &    &  0   & \frac{uv}{4t^2wK}\widetilde{C_1}   &\frac{K^2-v^2}{2Kvw}\widetilde{C_2}& -\frac{uw}{4t^2vK}\widetilde{B_1}  & \frac{-K^2+w^2}{2Kvw}\widetilde{B_2} \\
 & & & & & & & & \\
\widetilde{B_1}     &      &      &    &      & 0   &0 &\frac{v^2-w^2}{2vwt}\widetilde{A_1}+  &\frac{v^2-w^2}{2vwt}\widetilde{A_2} -\\
 & & & & & & &\frac{u(v^2-w^2)}{4t^2vwK}\widetilde{A_4}  & \frac{u(v^2-w^2)}{4t^2vwK}\widetilde{A_3} \\
 & & & & & & & & \\
 \widetilde{B_2}    &      &      &    &      &     & 0   &\frac{v^2-w^2}{2vwK}\widetilde{A_3} &\frac{v^2-w^2}{2vwK}\widetilde{A_4} \\
  & & & & & & & & \\
\widetilde{C_1}     &      &      &    &      &     &     & 0   &0 \\
 & & & & & & & & \\
\widetilde{C_2}   &      &      &    &      &     &     &     &  0  \\
\end{array}$$

\medskip

Infinitesimal isometries  are given by the vectors $X \in \m$ satisfying $$B(A_XY,Z)=-B(Y,A_XZ),$$ for any $Y, Z \in\m,$ where $A_XY=|X,Y]-\nabla_XY$. Since $\nabla_XY=U(X,Y)+\frac{1}{2}[X,Y]$, we have $A_XY=-U(X,Y)+\frac{1}{2}[X,Y]$. Thus $X \in\m$ is an infinitesimal isometry if
$$B([X,Y],Z)+B(Y,[X,Z])=0,$$
for any $Y,Z \in \m.$
\begin{proposition}
If $X \in \m$ is an infinitesimal isometry, thus
\begin{enumerate}
\item If $u=0$ and
\begin{enumerate}
\item If $t^2=v^2$, $X=c_1\widetilde{C_1}+c_2\widetilde{C_2},$
\item If $t^2=w^2$, $X=c_1\widetilde{B_1}+c_2\widetilde{B_2},$
\item If $v^2=w^2$, $X=a_1\widetilde{A_1}+a_2\widetilde{A_2}+a_3\widetilde{A_3}+a_4\widetilde{A_4}.$
\end{enumerate}
\item If $u \neq 0$ and
\begin{enumerate}
\item If $v^2=w^2$, $X=a_1\widetilde{A_1}+a_2\widetilde{A_2}+a_3\widetilde{A_3}+a_4\widetilde{A_4},$
\item If $v^2 \neq w^2$, $X=0$.
\end{enumerate}
\end{enumerate}
\end{proposition}
In particular, if $u=0$ and $t^2=v^2=w^2$, the bilinear form $B$ is $ad(\m)$-invariant.
\section{On the (non)naturally reductivity}
We consider on the \zz $\ SO(5)/SO(2)\times SO(2) \times SO(1)$ the $\Z^2_2$-symmetric Riemannian metric associated with the bilinear form
$$B=t^2(\alpha_1^2+ \alpha_2^2+\alpha_3^2+\alpha_4^2)+u(\alpha_1 \alpha_4-\alpha_2 \alpha_3 )+
v^2(\beta_1^2 +\beta_2^2)+w^2( \gamma_1^2+ \gamma_2^2)$$
with $tvw \ne 0$ and $u \in \left]-4t^2, 4t^2\right[$.
\begin{proposition}
The \zz $ \ SO(5)/SO(2)\times SO(2) \times SO(1)$ is naturally reductive if and only if
$$u=0 \ {\rm and} \ t^2=v^2=w^2.
$$
\end{proposition}
\noindent \pf Indeed, naturally reductivity means that
$$B(X,[Z,Y]_{\m})+B([Z,X]_{\m},Y)=0,$$
that is, $B(U(X,Y),Z)=0$ for any $X,Y,Z \in \m$.  From the expression of $U$ we deduce $u=0$ and $t^2=v^2=w^2=K^2.$ But $u=0$ implies $K^2=t^2.$ $\clubsuit$
\bigskip

Now, we assume that we have not the  naturally reductivity property. In this case, we can study if such a Riemannian space is a d'Atri space, that is, the geodesic symmetries preserve the volume. Let us recall that naturally reductive homogeneous spaces are d'Atri spaces, but these two notions are not equivalent. We know some examples of d'Atri spaces which are not naturally reductive. The condition of being a d'Atri space is often difficult to compute because it is necessary to know the equations of geodesics. But the D'Atri definition is equivalent to the Ledger (infinite) system whose the first equation writes
$$L(X,Y,Z)= (\nabla_{X}\rho)(Y, Z) + (\nabla_{Y}\rho)(Z,X) + (\nabla_{Z}\rho)(X,Y)= 0,$$
for any $X,Y,Z \in \m$, 
 where $\rho$ is the Ricci tensor, that is, the trace of the linear operator $V \rightarrow R(V,X)Y$. In the orthonormal basis previously defined the Ricci tensor is the symmetric matrix
$$\left(
        \begin{array}{cccccccc}
          \rho_{11} & 0 & 0 & \rho_{14} & 0 & 0 & 0 & 0 \\
          0 & \rho_{22} & \rho_{23} & 0 & 0 & 0 & 0 & 0 \\
          0 & \rho_{23} & \rho_{33} & 0 & 0 & 0 & 0 & 0\\
          \rho_{14} & 0 & 0 & \rho_{44} & 0 & 0 & 0 & 0 \\
          0 & 0 & 0 & 0 & \rho_{55} & 0 & 0 & 0 \\
          0 & 0 & 0 & 0 & 0 & \rho_{66} & 0 & 0 \\
          0 & 0& 0 &0 & 0 & 0 & \rho_{77} & 0 \\
          0 & 0 & 0 & 0 & 0 & 0 & 0 & \rho_{88} \\
        \end{array}
      \right)
$$
with
$$\begin{array}{l}
\medskip
    \rho_{11}=\rho_{22} =\frac{4t^4+u^2-4(v^4-6v^2w^2 + w^4)}{8t^2v^2w^2}, \\
    \medskip
    \rho_{14}=- \rho_{23}=\frac{u(K^2t^2 + v^4 - 6v^2w^2 + w^4)}{4 Kt^3v^2w^2},\\
    \medskip
    \rho_{33}= \rho_{44}=\frac{(4t^4-u^2)^2-4(v^4-6v^2w^2+w^4)(u^2+4t^4)}{8t^2(4t^4-u^2)v^2w^2},\\
    \medskip
    \rho_{55}= \rho_{66}=\frac{-4t^6+12t^4w^2+t^2(u^2+4v^4-4w^4)-3u^2w^2}{(4t^4-u^2)v^2w^2}, \\
    \medskip
    \rho_{77}= \rho_{88}=\frac{-4t^6+12t^4v^2+t^2(u^2-4v^4+4w^4)-3u^2v^2}{(4t^4-u^2)v^2w^2}.
  \end{array}
  $$
The non trivial Ledger equations (for the first condition) are
$$
\left\{
\begin{array}{l}
L(\widetilde{A_1},\widetilde{B_1},\widetilde{C_1})=L(\widetilde{A_2},\widetilde{B_1},\widetilde{C_2}),\\
L(\widetilde{A_1},\widetilde{B_2},\widetilde{C_2})=L(\widetilde{A_2},\widetilde{B_2},\widetilde{C_1}),\\
L(\widetilde{A_3},\widetilde{B_1},\widetilde{C_2})=L(\widetilde{A_4},\widetilde{B_1},\widetilde{C_1}),\\
L(\widetilde{A_3},\widetilde{B_2},\widetilde{C_1})=L(\widetilde{A_4},\widetilde{B_2},\widetilde{C_2}).
\end{array}
      \right.
$$
Thus we obtain
$$(*) \ \
\left\{
\begin{array}{l}
\medskip
(v^2-w^2) \rho_{11}+ (w^2-t^2) \rho_{55}+(t^2-v^2) \rho_{77}+\frac{u(w^2-v^2)}{2tK} \rho_{14}=0,\\
\medskip
-\frac{u}{2t} \rho_{55}+\frac{u}{2t} \rho_{77}+\frac{v^2-w^2}{K} \rho_{14}=0,\\
\medskip
\frac{u(v^2-w^2)}{2tvwK} \rho_{33}+\frac{uw}{2tvK} \rho_{55}-\frac{v^2-w^2}{vw} \rho_{14}-\frac{ u v}{2 t w K} \rho_{77}=0,\\
\medskip
(v^2-w^2) \rho_{33}+ (w^2-K^2) \rho_{55}+(K^2-v^2) \rho_{77}=0.
\end{array}
      \right.
$$
\begin{itemize}
\item If $u=0$, thus $ \rho_{11}=\rho_{33}$, $\rho_{14}=0$ and $(*)$ is equivalent to
$$(v^2-w^2) \rho_{33}+ (w^2-t^2) \rho_{55}+(t^2-v^2) \rho_{77}=0,$$
that is,
$$
\begin{array}{l}
(v^2-w^2)(9t^4+v^4+10v^2w^2+w^4-10t^2v^2-10t^2w^2)=0
\end{array}
$$
We obtain
$$\left\{
\begin{array}{l}
v^2=w^2 \\
 {\rm or} \
9t^4+v^4+10v^2w^2+w^4-10t^2v^2-10t^2w^2=0.
\end{array}
\right.$$
Let us study the second equation. For this, since $t \neq 0$, we can consider the change of variables
$$V=\frac{v^2}{t^2}, \ \ W=\frac{w^2}{t^2}.$$
The equation becomes
$$9  - 10(V+W) + (V+W)^2 +8VW= 0.$$
Now we put $S=V+W,$ $P=VW$ and we obtain
$$P=\frac{-S^2+10S-9}{8}.$$
Since $V$ and $W$ are strictly positive, $P > 0$, which implies that $-S^2+10S-9 >0$, that is, $S \in \left]1,9\right[.$
With these conditions,  $V$ and $W$ are the roots of $X^2-SX+P=0$. In fact they always exist since
$X^2-SX+P=X^2-SX+\frac{-S^2+10S-9}{8}$ has as discriminant
$$\Delta = S^2-4P=S^2-\frac{-S^2+10S-9}{2}=\frac{3S^2-10S+9}{2}$$
which is always positive. The roots of  $X^2-SX+P$ are
$$\left\{
\begin{array}{l}
X_1=\displaystyle\frac{S-\sqrt{\frac{3S^2-10S+9}{2}}}{2},\\
\medskip
X_2= \displaystyle\frac{S+\sqrt{\frac{3S^2-10S+9}{2}}}{2}.\\
\end{array}
\right.
$$
Since $S \in  \left]1,9\right[$, $P >0$ and $X_1 >0, X_2 >0$. We obtain
$$\left(\frac{v^2}{t^2}, \frac{w^2}{t^2} \right)=(X_1,X_2) \ \  {\rm{or}} \ (X_2,X_1).$$
\begin{proposition} Assume that the Riemannian $\Z^2_2$-symmetric metric on $SO(5)/SO(2) \times SO(2) \times SO(1)$ associated with a bilinear form $B$ with $u=0$ satisfies the first Ledger condition. Then
    $B$ is one of the following bilinear form
    $$
    \begin{array}{lll}
  B_1&=&t^2(\alpha_1^2+ \alpha_2^2+\alpha_3^2+\alpha_4^2)+v^2(\beta_1^2 +\beta_2^2+ \gamma_1^2+ \gamma_2^2),\\
B_2(S)&=&t^2(\alpha_1^2+ \alpha_2^2+\alpha_3^2+\alpha_4^2+\displaystyle\frac{S-\sqrt{\frac{3S^2-10S+9}{2}}}{2}(\beta_1^2 +\beta_2^2)\\
&& +\displaystyle\frac{S+\sqrt{\frac{3S^2-10S+9}{2}}}{2}( \gamma_1^2+ \gamma_2^2)),\\
B_3(S)&=&t^2(\alpha_1^2+ \alpha_2^2+\alpha_3^2+\alpha_4^2+\displaystyle\frac{S+\sqrt{\frac{3S^2-10S+9}{2}}}{2}(\beta_1^2 +\beta_2^2)\\
&& +\displaystyle\frac{S-\sqrt{\frac{3S^2-10S+9}{2}}}{2}( \gamma_1^2+ \gamma_2^2)),
   \end{array}
   $$
       with $S \in \left]1,9\right[.$  The metrics associated with $B_2(S)$ and $B_3(S)$ are not naturally reductive and the matrix associated with $B_1$ is not naturally reductive as soon as $t^2 \neq v^2$.
    \end{proposition}

\item Assume that $u \neq 0.$ The first Ledger condition writes

$$ \left\{
\begin{array}{l}
-t^2(v^2 - w^2)(36t^6 - 40t^4(v^2 + w^2) +6u^2(v^2 + w^2) \\
       \ \ \ \    +t^2(-9 u^2 + 4(v^4 + 10 v^2 w^2 + w^4)))=0, \\

t^2u(v^2 - w^2)(28 t^4 - 7u^2 - 16t^2(v^2 + w^2) +
        4(v^4 - 6 v^2w^2 + w^4))=0,  \\

 (v^2 - w^2)(144 t^8 - 160 t^6(v^2 + w^2) +
        40t^2u^2(v^2 + w^2) -\\
       \ \ \ \  16t^4(4 u^2 - v^4 - 10 v^2w^2 - w^4) +
        u^2(7u^2 - 4(v^4 - 6 v^2w^2 + w^4)))=0,  \\

                    t^2u(v^2 - w^2)
    (4 t^6 - 12 t^4(v^2 + w^2) + 3 u^2(v^2 + w^2)-
        t^2(u^2 - 32 v^2w^2))=0.
 \end{array}
\right.$$
We obtain $v^2=w^2$ or
$$ \left\{
\begin{array}{l}
36t^6 - 40t^4(v^2 + w^2) + 6u^2(v^2 + w^2) +
        t^2(-9 u^2 + 4(v^4 + 10 v^2 w^2 + w^4))=0, \\

28 t^4 - 7u^2 - 16t^2(v^2 + w^2) +
        4(v^4 - 6 v^2w^2 + w^4))=0,  \\

144 t^8 - 160 t^6(v^2 + w^2) +
        40t^2u^2(v^2 + w^2) -\\
       \ \ \ \  16t^4(4 u^2 - v^4 - 10 v^2w^2 - w^4) +
        u^2(7u^2 - 4(v^4 - 6 v^2w^2 + w^4))=0,  \\

    4 t^6 - 12 t^4(v^2 + w^2) + 3 u^2(v^2 + w^2)-
        t^2(u^2 - 32 v^2w^2)=0.
 \end{array}
\right.$$
\end{itemize}

The change of variables: $U=\frac{u}{t^2}, \ V=\frac{v^2}{t^2}, \ W=\frac{w^2}{t^2},\ S=V+W$ and $P=VW$
shows that the first three equations of the previous system are equivalent. So the system reduces to
the two equations
\begin{equation}
64P-24PS+4S-13S^2+3S^3=0, \label{Eq1}
\end{equation}
\begin{equation}
7U^2=28-16S+4(S^2-8P)\label{Eq2}.
\end{equation}

Equation \ref{Eq1} gives
$$P=\displaystyle\frac{-4S+13S^2-3S^3}{8(8-3S)}=\frac{-S(S-4)(3S-1)}{8(8-3S)},$$
because $S=\frac{8}{3}$ is not a solution. Thus $V$ and $W$ are roots of $X^2-SX+P=0$ if there are two real solutions of this equation.
But
$$X^2-SX+P=X^2-SX+\frac{-S(S-4)(3S-1)}{8(8-3S)}$$ and
its discriminant is $$\Delta= S^2-4P= S^2+\frac{S(S-4)(3S-1)}{2(8-3S)}=\frac{S(-3S^2+3S+4)}{2(8-3S)}.$$ We obtain a condition on $S$
to the existence of $V$ and $W$ which is
$$S \in \left]0, \frac{3+\sqrt{57}}{6}\right[\ \bigcup \ \left]\frac{8}{3},+\infty\right[.$$ Moreover we want $V$ and $W$ to be positive solutions so
$P>0$ and $S>0$; we have to take $S \in \left]\frac{1}{3}, \frac{8}{3}\right[\ \bigcup \ \left]4, +\infty\right[.$
Finally
$$S \in \left]\frac{1}{3}, \frac{3+\sqrt{57}}{6}\right[ \ \bigcup \ \left]4,+\infty\right[.$$

Then Equation \ref{Eq2} gives conditions on $S$ to the existence of $U$. In fact Equation \ref{Eq2} is equivalent to
$$U^2=4\frac{8-7S+S^2}{8-3S}.$$
Then $S \in \left]0,\frac{7-\sqrt{17}}{2}\right[ \bigcup \left] \frac{8}{3},\frac{7+\sqrt{17}}{2} \right[$. But we have also that $u \in \left]-4t^2 , 4t^2 \right[$
so $U^2<16$ which reduces to $S^2 +5S-24<0$ and $S\in  \left]0,\frac{7-\sqrt{17}}{2}\right[  \bigcup \left] \frac{8}{3},3 \right[$.

With the conditions coming from \ref{Eq1} we finally have that  $S \in \left]\frac{1}{3},\frac{7-\sqrt{17}}{2}\right[$

\begin{proposition}
Assume that the $\Z^2_2$-symmetric metric on $SO(5)/SO(2) \times SO(2) \times SO(1)$ associated with a bilinear form $B$ with $u\neq 0$ satisfies the first Ledger condition. Then $S$ belongs to $\left]\frac{1}{3},\frac{7-\sqrt{17}}{2}\right[$ and
    $B=t^2(\alpha_1^2+ \alpha_2^2+\alpha_3^2+\alpha_4^2)+u(\alpha_1\alpha_4 -\alpha_2\alpha_4)+v^2(\beta_1^2 +\beta_2^2)+w^2( \gamma_1^2+ \gamma_2^2),$ with
    \begin{enumerate}
    \item  $$v^2=\frac{-16 S + 6 S^2 -\sqrt{2S(32 + 12 S - 33 S^2 + 9S^3)}}{-32 + 12 S}t^2$$ and  $$ w^2=\frac{-3S(S-4)(S-\frac{1}{3})}{8(8-3S)V}t^2$$
    \item $$w^2=\frac{-16 S + 6 S^2 -\sqrt{2S(32 + 12 S - 33 S^2 + 9S^3)}}{-32 + 12 S}t^2 $$ and  $$ v^2=\frac{-3S(S-4)(S-\frac{1}{3})}{8(8-3S)V}t^2$$
      \end{enumerate}
    and in each case $u^2=4\frac{8-7S+S^2}{8-3S}t^4.$
    None of these metrics is  naturally reductive.
\end{proposition}

\end{document}